\documentclass[11pt]{amsart}

\usepackage[margin=1in]{geometry}
\usepackage{amsmath,amssymb,amsfonts,amsthm,mathtools}
\usepackage{mathrsfs}
\usepackage{enumitem}
\usepackage{hyperref}

\hypersetup{
  colorlinks=true,
  linkcolor=blue,
  citecolor=green,
  urlcolor=blue
}

\newtheorem{thm}{Theorem}[section]
\newtheorem{lem}[thm]{Lemma}
\newtheorem{prop}[thm]{Proposition}

\newtheorem{fact}[thm]{Fact}

\theoremstyle{definition}
\newtheorem{defn}[thm]{Definition}
\newtheorem{rmk}[thm]{Remark}
\newtheorem{ques}[thm]{Question} 

\DeclareMathOperator{\GL}{GL}

\DeclareMathOperator{\UT}{UT}
\DeclareMathOperator{\diag}{diag}

\newcommand{\NN}{\mathbb{N}}
\newcommand{\ZZ}{\mathbb{Z}}
\newcommand{\QQ}{\mathbb{Q}}
\newcommand{\RR}{\mathbb{R}}
\newcommand{\CC}{\mathbb{C}}

\begin{document}

\title{A faithful representation and profinite completion}

\author{Wonyong Jang}
\address{Department of Mathematical Sciences, Seoul National University, 1 Gwanak-ro, Gwanak-gu, Seoul 08826, Republic of Korea}
\email{jangwy@snu.ac.kr}

\begin{abstract}
We exhibit two finitely generated residually finite groups $G_1$ and $G_2$ with isomorphic profinite completions $\widehat{G_1} \cong \widehat{G_2}$ such that $G_1$ has a faithful representation into $\GL(d,\CC)$ but $G_2$ is not linear.
This answers a question of Andrei Jaikin-Zapirain, Alexander Lubotzky, and Martin Bridson negatively.
Moreover, $G_1$ and $G_2$ form a Grothendieck pair such that every faithful finite-dimensional complex representation of $G_1$ extends uniquely to a nonfaithful representation of $G_2$.
We also present two finitely generated residually finite groups $H_1$ and $H_2$ with isomorphic profinite completions $\widehat{H_1} \cong \widehat{H_2}$ satisfying the following.
For every field $k$ with characteristic $p \geq 5$ and transcendence degree $r \geq 2$,
there exists $d>0$ such that $H_1$ has a faithful representation into $\GL(d,k)$. 
However, $H_2$ admits no faithful finite-dimensional representation over any field of positive characteristic.
\end{abstract}

\maketitle

\vspace{0.6cm}
\noindent \textbf{AMS Classification numbers (2020)} Primary: 20H20; Secondary: 20E18, 20E26

\vspace{0.6cm}

\section{Introduction}\label{sec:intro}

 A group $G$ is called \emph{linear} if $G$ has a faithful representation into $\GL(d,k)$ for some $d>0$ and some field $k$.
 Finitely generated linear groups share many important properties; for example, they are residually finite, and if the base field has characteristic $0$, then they are virtually torsion-free.
 Moreover, they satisfy the Tits alternative; a finitely generated linear group either contains a rank $2$ free group or is virtually solvable.
 This celebrated result implies that the von Neumann conjecture --- $G$ is non-amenable if and only if $G$ contains a rank $2$ free group --- is true for linear groups.

 From the above perspective, the attempt to determine whether a given group is linear continues, and many results have been established.
 For instance, the surface group $\pi_1(S_g)$ and braid groups are linear.
 Recently, it was also proved that the fundamental groups of compact 3-manifolds \cite{wu20263} and graph manifold groups \cite{wise2026linearity} are linear.
 However, there exist hyperbolic groups which are not linear (\cite[Section 8]{MR2175160}, \cite{canary2019new}, and \cite{tholozan2025residually}).
 Also, Thompson's groups $F,T,V$ and Grigorchuk's group are not linear.

 For each group property $\mathcal{P}$, one can ask whether the property $\mathcal{P}$ is determined by looking at finite quotients.
 This question is motivated by the \emph{profinite rigidity problem} concerning the classification of finitely generated residually finite groups by their finite quotients.
 Since the profinite completion $\widehat{G}$ encodes all finite quotients of $G$, we can consider the ``profinite invariant'' group property as follows.
 A group property $\mathcal{P}$ is called a \emph{profinite property} if for any two finitely generated residually finite groups $\Gamma_1$ and $\Gamma_2$ with $\widehat{\Gamma}_1 \cong \widehat{\Gamma}_2$, the group $\Gamma_1$ has $\mathcal{P}$ if and only if $\Gamma_2$ has $\mathcal{P}$.

While some properties, such as satisfying a nontrivial law \cite{bridson2026chasing} and being polycyclic \cite{sabbagh1991polycyclic}, are profinite properties, it is well known that many geometric and asymptotic properties are not.
Notable examples of properties that fail to be profinite include amenability \cite{kionke2023amenability}, 
Kazhdan’s property (T) \cite{aka2012profinite}, 
being conjugacy separable \cite{cotton2012conjugacy}, 
various finiteness properties \cite{lubotzky2014finiteness}, 
Serre’s property FA \cite{cheetham2024property}, 
second bounded cohomology \cite{echtler2024bounded}, 
being torsion-free \cite{bridson2024profinite}, 
two orderabilities (\cite{bridson2024profinite} and \cite{jang2026bi}), and
being co-Hopfian \cite{baik2026co}.
 Recently, it has also been proved that stable commutator length and quasimorphisms are not profinite invariants \cite{fournier2026profinite}.

 Although there are examples of profinite and non-profinite properties, 
 for many interesting properties, it remains unknown whether they are profinite.
 In particular, the case of linearity is of particular interest.
 The following remarkable result was proved in \cite{grothendieck1970representations} and is also discussed in \cite{jaikin2025some}.

 \begin{thm}[Théorème 1.2 (b) in \cite{grothendieck1970representations} and Theorem 1.1 (2) in \cite{jaikin2025some}] \label{thm:Grothendieck}
 Let $H \hookrightarrow G$ be a Grothendieck pair.
 For every field $k$ and every $d \geq 1$, each representation $\rho:H \to \GL(d,k)$ admits a unique extension $\widetilde{\rho}:G \to \GL(d,k)$.
 \end{thm}
 Recall that a pair of finitely generated residually finite groups $H \hookrightarrow G$ is called a \emph{Grothendieck pair} if the inclusion induces $\widehat{H} \cong \widehat{G}$.
 Here, $\widehat{G}$ is the profinite completion of $G$, which we recall in Section~\ref{subsec:profinite}.
 It is then natural to ask whether the unique extension of a faithful representation must remain faithful.
 A related question can be found in \cite{jaikin2025some} and \cite[(11) in Section 5.2]{bridson2026chasing}.

\begin{ques} \label{conj:Lpp}
 Is linearity a profinite property?
\end{ques}

 Our first main theorem gives a negative answer to this question.

 \begin{thm}[Theorem \ref{thm:main1}] \label{thm:intro-main1}
Linearity is not a profinite property.
In other words, there exist finitely generated residually finite groups $G_1$ and $G_2$ such that:
\begin{enumerate}[label=(\arabic*)]
  \item $\widehat{G_1} \cong \widehat{G_2}$;
  \item $G_1$ has a faithful representation into $\GL(d,\CC)$ for some $d>0$;
  \item $G_2$ does not have a faithful representation into $\GL(d,k)$ for any field $k$ and any $d>0$.
\end{enumerate}
\end{thm}

 In our construction, $G_1$ and $G_2$ are indeed a Grothendieck pair.
 By Grothendieck’s theorem (Theorem~\ref{thm:Grothendieck}), 
 every faithful finite-dimensional complex representation of $G_1$ extends uniquely to $G_2$.
 Since $G_2$ is not linear, every such extension has a nontrivial kernel. 
 Thus, faithfulness need not be preserved under the unique extension associated with a Grothendieck pair.
 The author is not aware of any earlier examples of this phenomenon.
 Theorem~\ref{thm:intro-main1} also implies that admitting a faithful representation into $\GL(d,k)$ with $\textnormal{char}(k)=0$ is not a profinite property.
 
 Next, we focus on a faithful representation into $\GL(d,k)$ for a fixed field with  $\textnormal{char}(k)=p \geq 5$.
 Our second main theorem shows that having a faithful representation into $\GL(d,k)$ is not a profinite property for every field $k$ of characteristic $p \geq 5$ and transcendence degree $r \geq 2$.

\begin{thm}[Theorem \ref{thm:main2}] \label{thm:intro-main2}
There exist two finitely generated residually finite groups $H_1$ and $H_2$ such that:
\begin{enumerate}[label=(\arabic*)]
  \item $\widehat{H_1} \cong \widehat{H_2}$;
  \item For every field $k$ with characteristic $p \geq 5$ and transcendence degree $r \geq 2$,
there exists $d>0$ such that $H_1$ has a faithful representation into $\GL(d,k)$;
  \item $H_2$ cannot admit a faithful representation into $\GL(d,L)$ for any $d>0$ and any field $L$ with characteristic $p > 0$.
\end{enumerate}
\end{thm}

 In our construction, $H_1<H_2$ is again a Grothendieck pair. 
 Consequently, over every field $F$ of characteristic $p \geq 5$ and transcendence degree $r \geq 2$, $H_1$ has a faithful finite-dimensional representation whose unique extension to $H_2$ is not faithful.

 The constructions in Theorem~\ref{thm:intro-main1} and Theorem~\ref{thm:intro-main2} are independent, but both constructions rely on an ``iterated'' fiber product.
 Recall that for given two surjective homomorphisms $f_1:A \to Q_1$ and $f_2:B \to Q_1$, the corresponding fiber product is $P_1:=\{ (a,b) \in A \times B : f_1(a)=f_2(b) \}$.
 Let $\pi_1:P_1 \to B$ be the projection onto the second factor.
 For a surjective homomorphism $\pi_2:C \to B$, we can define the fiber product $P_2$ of $\pi_1$ and $\pi_2$.
 The key point is that $P_2$ is actually a subgroup of $A \times C$ since
 \begin{align*}
 P_2 &= \{ (x,c) \in P_1 \times C : \pi_1(x)=\pi_2(c) \} \\
 & \cong \{ (a,b,c) \in A \times B \times C : \pi_1((a,b))=\pi_2(c), f_1(a)=f_2(b) \} \\
 & \cong \{ (a,b,c) \in A \times B \times C : b=\pi_2(c), f_1(a)=f_2(b) \} \\
 & \cong \{ (a,c) \in A \times C : f_1(a)=f_2(\pi_2(c)) \}.
 \end{align*}
 Under suitable assumptions, we obtain $\widehat{P_2} \cong \widehat{P_1 \times C}$.
We then choose linear groups $A$ and $C$ so that $P_1$ is not linear.
Under the above assumptions, the groups $P_2$ and $P_1 \times C$ provide a counterexample to Question~\ref{conj:Lpp}.
 We also use ``center-preserving'' embeddings introduced in \cite{houcine2007embeddings} and \cite{palmer2025embedding} to establish Theorem~\ref{thm:intro-main1}.
 The main obstruction to linearity is that the center contains $\oplus_\ZZ \ZZ/2 \ZZ$.
 This rules out linearity in characteristic $p \neq 2$; a separate argument handles characteristic $2$ (see Proposition~\ref{prop:Not_linear}).
 This is the reason why we use such special embeddings.
 To prove Theorem~\ref{thm:intro-main2}, we use the universal central extension $E$ of Thompson's group $T$.
 This is because there exists a surjection from a surface group $\pi_1(S_g) \to T$, and $\pi_1(S_g)$ has a faithful representation into $\GL(d,k)$ for some $d>0$ and any field $k$ with characteristic $p \geq 5$ and transcendence degree $r \geq 2$ \cite[Theorem 5]{fisher2018character}.
 To rule out faithful representations in positive characteristic, 
 we use the fact that $H_2(T;\ZZ) \cong \ZZ^2$ provides an infinite-order central element with the desired property
 (see also Lemma~\ref{Lem:Obstruction1}).
 Since these two constructions of Theorem~\ref{thm:main1} and Theorem~\ref{thm:main2} are independent, we reuse the notation in Section~\ref{sec:main1} and Section~\ref{sec:main2}.

The paper is organized as follows.
In Section~\ref{sec:prelim}, we briefly recall the profinite completion and the fiber product.
In addition, we introduce the main ingredients of the main result, the labelled Thompson group and the universal central extension of Thompson's group $T$.
We prove the first main result, Theorem~\ref{thm:main1}, in Section~\ref{sec:main1} and the second main theorem, Theorem~\ref{thm:main2}, in Section~\ref{sec:main2}.

\vspace{0.4cm}

\noindent \textbf{Acknowledgment.}
 We would like to thank Hyungryul Baik, Andrei Jaikin-Zapirain, Alex Lubotzky, and Alan Reid for their interest in this work and encouraging comments.
 The author also thanks Francesco Fournier-Facio and Sam Hughes for valuable comments that improved the manuscript.
 This research was supported by BK21 SNU Mathematical Sciences Division.

\vspace{0.4cm}

\noindent \textbf{AI disclosure.}
The idea of using an ``iterated'' fiber product and the strategy to find a counterexample to Question~\ref{conj:Lpp} originated with the author.
To find some groups with the desired properties, the author used ChatGPT 6.1 Sol.
The author also used ChatGPT for finding references and grammar checking.
The author takes full responsibility for this final version.

\section{Preliminaries}\label{sec:prelim}

\subsection{Profinite completion and fiber product} \label{subsec:profinite}
 One of the main concepts in this paper is that of profinite properties.
  To introduce this notion precisely, we begin with the definition of the profinite completion.

\begin{defn}
 For a group $G$, let $\mathcal{N}$ be the set of all normal subgroups of finite index in $G$.
Consider the following inverse limit
\[
\widehat{G}=\varprojlim_{N\in\mathcal{N}} G/N.
\]
The group $\widehat{G}$ is called the \emph{profinite completion} of $G$.
\end{defn}
 
\begin{rmk}
 There is a canonical homomorphism from $G$ to its profinite completion $\widehat{G}$ given by \[ i : g \mapsto \{ gN \}_{N \in \mathcal{N}}. \]
This homomorphism is not injective in general. Indeed, the homomorphism $i$ is injective if and only if $G$ is residually finite.
\end{rmk}

Let $\mathcal{C}(G)$ be the set of isomorphism classes of all finite quotients of $G$.
For finitely generated groups $G$ and $H$, $\widehat{G} \cong \widehat{H}$ if and only if $\mathcal{C}(G) = \mathcal{C}(H)$ (see \cite{dixon1982profinite} and \cite{nikolov2007finitely}).
In other words, $G$ and $H$ have isomorphic profinite completions if and only if they have the same finite quotients up to isomorphism.
Now we can use the concept of a profinite completion to introduce the profinite property.

\begin{defn}
 Let $\mathcal{P}$ be a group property of finitely generated residually finite groups. We say that $\mathcal{P}$ is a \emph{profinite property} if for any two finitely generated residually finite groups $G$ and $H$ with $\widehat{G} \cong \widehat{H}$, $G$ satisfies $\mathcal{P}$ if and only if $H$ satisfies $\mathcal{P}$.
\end{defn}

Namely, if a group property $\mathcal{P}$ is a profinite property, then it can be detected by profinite completion, and equivalently, by taking finite quotients.
For more details on profinite groups, profinite rigidity, and profinite properties, 
we refer to \cite{nikolov2007finitely}, \cite{MR2599132}, \cite{reid2018profinite}, and \cite{bridson2026chasing}.

We close this subsection by recalling the fiber product.
We also record a well-known criterion ensuring that the profinite completions of the fiber product and the direct product are isomorphic.

\begin{defn}\label{def:fiberproduct}
 Let $p_1 \colon G_1 \to Q$, $p_2\colon G_2 \to Q$ be two epimorphisms between groups. The \emph{fiber product} of $p_1$ and $p_2$ is the subgroup given by
 $$ P\coloneqq \{ (g_1,g_2) : p_1(g_1)=p_2(g_2) \} < G_1 \times G_2 . $$
\end{defn}

 Proposition~\ref{Fiber_profinite_completion} is a key ingredient in our constructions and can be seen in \cite[Proposition 1.2]{bridson2024profinite}.
 The original statement and its proof can be found in \cite[Theorem 2.2]{MR3531665}, and this type of proposition was first suggested by Platonov and Tavgen in \cite{platonov1986grothendieck}.
 It has also played an important role in providing an example of Grothendieck pairs. See \cite{MR2119723}.
 First, we record the following well-known fact.

\begin{lem}[Lemma 1.1 in \cite{bridson2024profinite}] \label{lem:P-f.g.}
Let $p_1 : G_1 \to Q$, $p_2 : G_2 \to Q$ be two epimorphisms of groups, and $P$ be the fiber product of $p_1$ and $p_2$.
If both $G_1$ and $G_2$ are finitely generated and $Q$ is finitely presented, then $P$ is finitely generated.
\end{lem}

 Under additional assumptions, we can obtain the isomorphism between their profinite completions.

\begin{prop}[Proposition 1.2 in \cite{bridson2024profinite}] \label{Fiber_profinite_completion}
Let $p_1 : G_1 \to Q$, $p_2 : G_2 \to Q$ be two epimorphisms of groups, and $P$ be the fiber product of $p_1$ and $p_2$. 
If $G_1, G_2$ are finitely generated, $Q$ is finitely presented, $\widehat{Q}=1$, and $H_2(Q,\ZZ)=0$, then the inclusion $P \hookrightarrow G_1 \times G_2 $ induces an isomorphism of profinite completions $\widehat{P} \cong \widehat{G_1} \times \widehat{G_2}.$
\end{prop}

\subsection{Labelled Thompson groups} \label{subsec:Labelled Thompson groups}

In this subsection, we briefly recall the labelled Thompson group.
For a discrete group $G$, we can construct the labelled Thompson group $V(G)$.
This concept was first introduced in \cite{thompson1980embeddings} and has been well studied in \cite{wu2025embedding} and \cite{palmer2025embedding}.
In particular, there are several embeddings $\iota_0, \iota_{\emptyset}:G \to V(G)$.
For further details, we refer to these articles.
We record the following interesting results.
Throughout this article, we fix the embedding $\iota_{\emptyset}$ introduced in \cite{palmer2025embedding}.
So, by $G<V(G)$, we mean $\iota_\emptyset(G)<V(G)$.
For a subgroup $H<G$, we will write $H<V(G)$ instead of $\iota_\emptyset(H)<V(G)$.
 \begin{prop} \label{prop:G-V(G)}
 Let $G$ be a discrete group.
 Then the group $V(G)$ has the following properties.
 \begin{itemize}
    \item if $G$ is finitely presented, then $V(G)$ is also finitely presented \cite[Corollary 3.19]{wu2025embedding},
    \item $\widehat{V(G)}=1$ \cite[Corollary 2.3]{wu2025embedding},
    \item $H_2(V(G);\ZZ)=0$ \cite[Theorem 2.6]{palmer2025embedding}, and
    \item $Z(V(G)) \cong Z(G)$ \cite[Proposition 2.8]{palmer2025embedding}.
 \end{itemize}
 \end{prop}

\subsection{Thompson's group $T$ and its universal central extension} \label{subsec:UCET}

In this subsection, we aim to discuss the universal central extension of $T$.
Thompson’s groups $F,T,V$ were introduced by Richard Thompson in 1965.
We refer to \cite{cannon1996introductory} for a well-written introductory explanation of these groups.

Another main ingredient is the universal central extension of Thompson's group $T$.
We briefly recall the universal central extension and related results.
For any group $G$, an extension $$1 \to N \to E \to G \to 1$$ is called \emph{central} if $N$ is contained in the center $Z(E)$.
Let 
\begin{align*}
1 \to N_1 \to E_1 \overset{\pi_1}{\to} G \to 1 \\
1 \to N_2 \to E_2 \overset{\pi_2}{\to} G \to 1
\end{align*}
be two central extensions of $G$.
A homomorphism over $G$ from $E_1$ to $E_2$ is a homomorphism $f:E_1 \to E_2$ such that $\pi_2 \circ f = \pi_1$.
We say that $1 \to N \to E \to G \to 1$ is the \emph{universal central extension} if it is a central extension and for every central extension of $G$, $$ 1 \to N' \to E' \to G \to 1, $$
there exists a unique homomorphism $f:E \to E'$ over $G$.
It is clear that the universal central extension is unique, if it exists.
Also, the following is well known.

\begin{prop}[Theorem 6.9.5 in \cite{Weibel_1994}] \label{prop:AboutUCE}
 A group $G$ has the universal central extension if and only if $G$ is perfect.
 In this case, we have $$ 1 \to H_2(G;\ZZ) \to E \to G \to 1 $$
 where $E$ is the universal central extension of $G$.
\end{prop}
We refer to \cite[Section 6.9]{Weibel_1994} for readers who are interested in the universal central extension.
From Proposition~\ref{prop:AboutUCE} and \cite{heuer2021spectrum}, the universal central extension of $T$ has the following form.

\begin{fact}
 Let $T$ be Thompson's group and $E$ be its universal central extension.
 Then the corresponding central extension is given by
 $$ 1 \to \ZZ^2 \to E \overset{\pi}{\to} T \to 1. $$
\end{fact}

 Recall that $\ZZ^2$ comes from $H_2(T;\ZZ) \cong \ZZ^2$.
 We emphasize that we use the fact $H_2(T;\ZZ) \cong \ZZ^2$ (more precisely, $H_2(T;\ZZ)$ is torsion-free) to establish Theorem~\ref{thm:main2}.
 For readers who are interested in this extension of $T$, we refer to \cite[Section 5]{heuer2021spectrum}.
We record some properties of $E$ which will be used later.

\begin{prop}[Proposition 5.1 in \cite{heuer2021spectrum}] \label{prop:E-property}
 Let $E$ be the universal central extension of $T$. Then
 \begin{itemize}
    \item $E$ is finitely presented, and
    \item $H_1(E;\ZZ)=H_2(E;\ZZ)=0$.
 \end{itemize}
\end{prop}

 We also need the following elementary observation.

\begin{lem} \label{lem:E:trivial_completion}
The group $E$ does not have a proper finite index subgroup. In other words, $\widehat{E}=1$.
\end{lem}
\begin{proof}
 Suppose that $f:E \to F$ is a homomorphism where $F$ is finite.
 Let $C:= \ZZ^2 < Z(E)$, $I:=f(E)$ and $J:=f(C)$.
 Then $J < Z(I)$. From the homomorphism $g:E \to I \to I/J$, we get $g(C)=1$, so this induces a surjective homomorphism $h:E/C \to I/J$.
 Note that $E/C \cong T$, Thompson's group, so we have $h:T \to I/J$. 
 Since $\widehat{T}=1$ (recall that $T$ is an infinite simple group), $I/J$ is finite, and $h$ is surjective, we conclude that $I/J$ is trivial.
 This means that $I=J$ and $I$ is abelian.

 Since $E$ is perfect, so is $I=f(E)$.
 Thus, $I$ is abelian and perfect. This implies that $I$ is trivial.
 Hence, $f$ is trivial.
 Since $F$ is arbitrary, we complete the proof.
\end{proof}

 Combining Proposition~\ref{prop:E-property}, Lemma~\ref{lem:E:trivial_completion} allows us to apply Proposition~\ref{Fiber_profinite_completion} with $Q=E$.
 We notice that the idea of using Thompson groups to produce Grothendieck pairs first appears in \cite{corson2026higman}.
 Thompson’s group $V$ and the Higman–-Thompson groups $V_n$ are used in \cite{corson2026higman}, whereas the main ingredients in Section~\ref{sec:main2} are Thompson's group $T$ and its universal central extension $E$.
 In Section~\ref{sec:main1}, we use the labelled Thompson groups discussed in Section~\ref{subsec:Labelled Thompson groups}.

\section{Linearity and profinite completion} \label{sec:main1}
 We devote this section to proving the following main result.

\begin{thm} \label{thm:main1}
Linearity is not a profinite property.
In other words, there exist finitely generated residually finite groups $G_1$ and $G_2$ such that:
\begin{enumerate}[label=(\arabic*)]
  \item $\widehat{G_1} \cong \widehat{G_2}$;
  \item $G_1$ has a faithful representation into $\GL(d,\CC)$ for some $d>0$;
  \item $G_2$ does not have a faithful representation into $\GL(d,k)$ for any field $k$ and any $d>0$.
\end{enumerate}
\end{thm}

\subsection{The construction of two groups with isomorphic profinite completions}

 We start with a ring $R:=\ZZ[t,t^{-1}]$, and consider
 $$L:=\left< s:=
 \begin{bmatrix}
 t & 0 & 0 \\ 0 & 1 & 0 \\ 0 & 0 & t
 \end{bmatrix}, x:= \begin{bmatrix}
 1 & 1 & 0 \\ 0 & 1 & 0 \\ 0 & 0 & 1
 \end{bmatrix}, y:= 
 \begin{bmatrix}
 1 & 0 & 0 \\ 0 & 1 & 1 \\ 0 & 0 & 1
 \end{bmatrix}
 \right> < \GL(3,\QQ(t)).$$

 \begin{lem} \label{lem:L structure}
The group $L$ is isomorphic to $\UT(3,R) \rtimes \left< s \right>$ where $$\UT(3,R)=\left \{ 
\begin{bmatrix}
1 & a & c \\ 0 & 1 & b \\ 0 & 0 & 1
\end{bmatrix} : a,b,c \in R
\right \},$$
and the action is given by $$s.\begin{bmatrix}
1 & a & c \\ 0 & 1 & b \\ 0 & 0 & 1
\end{bmatrix} := s\begin{bmatrix}
1 & a & c \\ 0 & 1 & b \\ 0 & 0 & 1
\end{bmatrix}s^{-1}=\begin{bmatrix}
1 & ta & c \\ 0 & 1 & t^{-1}b \\ 0 & 0 & 1
\end{bmatrix}.$$
\end{lem}
\begin{proof}
 Clearly, $\left< s \right> < L$.
 From the action, $\UT(3,R) \cap L$ is a normal subgroup of $L$.
 We claim that $L$ contains $\UT(3,R)$.

 Recall the following elementary calculations:
 \begin{align*}
 \begin{bmatrix} 1 & f(t) & 0 \\ 0 & 1 & 0 \\ 0 & 0 & 1 
 \end{bmatrix}^{-1} \begin{bmatrix} 1 & 0 & 0 \\ 0 & 1 & g(t) \\ 0 & 0 & 1 
 \end{bmatrix}^{-1} \begin{bmatrix} 1 & f(t) & 0 \\ 0 & 1 & 0 \\ 0 & 0 & 1 
 \end{bmatrix} \begin{bmatrix} 1 & 0 & 0 \\ 0 & 1 & g(t) \\ 0 & 0 & 1 
 \end{bmatrix} &= \begin{bmatrix} 1 & 0 & f(t)g(t) \\ 0 & 1 & 0 \\ 0 & 0 & 1 
 \end{bmatrix} \\
 \begin{bmatrix} 1 & 0 & f(t) \\ 0 & 1 & 0 \\ 0 & 0 & 1 
 \end{bmatrix} \begin{bmatrix} 1 & 0 & 0 \\ 0 & 1 & g(t) \\ 0 & 0 & 1 
 \end{bmatrix} \begin{bmatrix} 1 & h(t) & 0 \\ 0 & 1 & 0 \\ 0 & 0 & 1 
 \end{bmatrix} &= \begin{bmatrix} 1 & h(t) & f(t) \\ 0 & 1 & g(t) \\ 0 & 0 & 1 
 \end{bmatrix}
 \end{align*}
 for $f(t),g(t),h(t) \in R$.
 Thus, to show that $L$ contains $\UT(3,R)$, we will show
 $$ \begin{bmatrix} 1 & t^n & 0 \\ 0 & 1 & 0 \\ 0 & 0 & 1 
 \end{bmatrix}, \begin{bmatrix} 1 & 0 & 0 \\ 0 & 1 & t^m \\ 0 & 0 & 1 
 \end{bmatrix} \in L$$
 for all $m,n \in \ZZ$.
 This fact is obtained from
 $$ s^nxs^{-n}=\begin{bmatrix} 1 & t^n & 0 \\ 0 & 1 & 0 \\ 0 & 0 & 1 
 \end{bmatrix}, \quad s^mys^{-m}=\begin{bmatrix} 1 & 0 & 0 \\ 0 & 1 & t^{-m} \\ 0 & 0 & 1 
 \end{bmatrix}.$$
 
 Thus, the subgroup generated by $s$ and $\UT(3,R)$ is contained in $L$.
 However, $s,x,y \in \UT(3,R) \rtimes \left< s \right>$, thus $L=\UT(3,R) \left< s \right>$.
 Also, note that $\left< s \right> \cap \UT(3,R)=\{ I_3 \}$.
 Therefore, $L$ is isomorphic to $\UT(3,R) \rtimes \left< s \right>$.
\end{proof}
 
 From Lemma~\ref{lem:L structure}, one can easily check that
 $$ Z(L)=\left \{  \begin{bmatrix}
 1 & 0 & z \\ 0 & 1 & 0 \\ 0 & 0 & 1
 \end{bmatrix} : z \in R \right \}. $$
 Thus, $$ 2Z(L):= \left \{  \begin{bmatrix}
 1 & 0 & 2z \\ 0 & 1 & 0 \\ 0 & 0 & 1
 \end{bmatrix} : z \in R \right \} $$ is a normal subgroup of $L$.
 Set 
 $$ H:=L/2Z(L). $$

\begin{lem} \label{lem:H-fgrp}
 The group $H$ is finitely generated and recursively presented.
\end{lem}
\begin{proof}
Obviously, $H$ is finitely generated.
Note that if a group $G$ has a solvable word problem, then $G$ is recursively presented \cite{rauzy2020computability}.
Hence, we will show that $H$ has a solvable word problem.

By Lemma~\ref{lem:L structure}, $L$ is isomorphic to $\UT(3,R) \rtimes \left< s \right>$.
This implies that every element $l$ of $L$ can be expressed as 
$$l=\left( \begin{bmatrix} 1 & h(t) & f(t) \\ 0 & 1 & g(t) \\ 0 & 0 & 1 
 \end{bmatrix},s^n \right ).$$
 For simplicity, write $l=(h(t),g(t),f(t),s^n)$.
 Then the multiplication on $L$ is given by

\begin{align} \label{eq:operationL}
 \begin{split}
&(h_1(t),g_1(t),f_1(t),s^{n_1}).(h_2(t),g_2(t),f_2(t),s^{n_2}) \\
=&(h_1(t)+t^{n_1}h_2(t),g_1(t)+t^{-n_1}g_2(t),f_1(t)+f_2(t)+h_1(t)t^{-n_1}g_2(t),s^{n_1+n_2})
  \end{split}
 \end{align}
 
 Since $H \cong L/2Z(L)$, every element $a \in H$ can be expressed as $a=(h(t),g(t),\overline{f(t)},s^n)$ where $\overline{f(t)} \in R/2R \cong \mathbb{F}_2[t,t^{-1}]$.
 From Formula~\eqref{eq:operationL}, the operation on $H$ is given by
 \begin{align} \label{eq:operationH}
 \begin{split}
  &(h_1(t),g_1(t),\overline{f_1(t)},s^{n_1}).(h_2(t),g_2(t),\overline{f_2(t)},s^{n_2}) \\ 
  =&(h_1(t)+t^{n_1}h_2(t),g_1(t)+t^{-n_1}g_2(t),\overline{f_1(t)}+\overline{f_2(t)}+\overline{h_1(t)t^{-n_1}g_2(t)},s^{n_1+n_2}).
 \end{split}
 \end{align}
 Here, the third coordinate is computed modulo $2R$.
 Choose an arbitrary word $w=w_1 \cdots w_\ell$ in $H$ where $w_i \in \{ s,x,y,s^{-1},x^{-1},y^{-1}\}$.
 Then each $w_i$ can be written as $(h_i(t),g_i(t),\overline{f_i(t)},s^{n_i})$.
 Using Formula~\eqref{eq:operationH}, we get a unique expression $w=(h(t),g(t),\overline{f(t)},s^{n})$.
 It is obvious that $w=1$ if and only if $h(t)=g(t)=\overline{f(t)}=n=0$.
 Thus the word problem of $H$ is solvable.
\end{proof}

 Combining with Lemma~\ref{lem:H-fgrp} and Theorem II in \cite{houcine2007embeddings}, there exists a finitely presented group $K$ such that $K$ contains $H$ as a subgroup, and $Z(K)=Z(H)$.
 Now consider the labelled Thompson group $M:=V(K)$ and the embedding $$\iota_\emptyset : K \hookrightarrow M.$$
 Due to Proposition~\ref{prop:G-V(G)}, $M$ is finitely presented, and both $\widehat{M}$ and $H_2(M;\ZZ)$ are trivial.
 In addition, $Z(M)=Z(K)$.
 For simplicity, let $$ C:=Z(M)=Z(K)=Z(H). $$
 Consider $$ A:=M/Z(M), \quad B:= H / Z(H). $$
 Since $H<M$, there exists an embedding from $B$ to $A$.
 Since $M$ is finitely generated, there exists a surjection $\sigma_1:F_{m_1} \to M$.
 Let $\Gamma := B*F_{m_1} $ and define $$r:\Gamma \to A$$
 where $r|B$ is the embedding from $B$ to $A$ and $r|F_{m_1}:F_{m_1} \twoheadrightarrow M \twoheadrightarrow M/Z(M) \cong A$.
 Hence, $r$ is surjective.

 Now we construct groups using an iterated fiber product.
 Let $$P_1 := \left \{ (\gamma,m) \in \Gamma \times M : r(\gamma)=mC \right \}$$ be the fiber product of $r:\Gamma \to A \cong M/C$ and the natural surjection $M \to M/C$.
 Let $q: P_1 \to M$ be the projection onto the second factor and choose a surjection $\sigma_2:F_{m_2} \to M$. 
 Consider the fiber product of $q$ and $\sigma_2$, $$P_2 := \left \{ (a,f) \in P_1 \times F_{m_2} : q(a)=\sigma_2(f) \right \}. $$
 Clearly, $F_{m_2}$ is finitely generated, and $M$ is finitely presented with $\widehat{M}=1$ and $H_2(M;\ZZ)=0$.
 To apply Proposition~\ref{Fiber_profinite_completion}, we only need to show that $P_1$ is finitely generated.
 To do this, we need the following lemma.

 \begin{lem} \label{lem:A amalgam}
 The group $P_1$ is isomorphic to $H*_C(C \times F_{m_1} )$.
 \end{lem}
 \begin{proof}
 Define two maps $$f:H \to P_1, \quad g: C \times F_{m_1} \to P_1 $$
 by $h \mapsto (hC,h)$ and $(c,w) \mapsto (w,c \sigma_1(w))$ (here, $\sigma_1:F_{m_1} \to M$), respectively.
 We will verify that both $f$ and $g$ are well-defined injective group homomorphisms.

 \begin{enumerate}
    \item First, consider $f$. \\
    $hC \in B$ so $r(hC)=hC$. This means $(hC,h) \in P_1$.
    Since the projection of $f$ onto the second factor is the identity map on $H$, $f$ is injective. By construction, $f$ is a homomorphism.
    \item Next consider $g$. \\ 
    Recall that $r(w)$ is defined as $\sigma_1(w)C$ (since $w \in F_{m_1}$), so this means that $(w,c \sigma_1(w))$ is an element in $P_1$.
    Note that $g$ is a homomorphism since $C = Z(M)$.
    Also, $g$ is injective since when $(w,c \sigma_1(w))=(1,1)$, we have $w=1$ and $c=1$.
 \end{enumerate}

 Let $p: P_1 \to \Gamma = B*F_{m_1}$ be the projection onto the first factor.
 We claim that
 $$ p^{-1}(B)=f(H), \quad p^{-1}(F_{m_1}) = g(C \times F_{m_1}). $$
 Since $hC \in B$ for $h \in H$ and $w \in F_{m_1}$, we have $ p^{-1}(B) \supset f(H) $ and $ p^{-1}(F_{m_1}) \supset g(C \times F_{m_1}). $
 For the opposite inclusions, choose $(b,m) \in p^{-1}(B) \subset P_1 \subset \Gamma \times M$.
 Since $B=H/C$, $b=hC$ for some $h \in H$.
 Since $(b,m) \in P_1$, by definition of the fiber product, we get $r(b)=mC$. Since $b \in B$, $r(b)=b$ so $b=mC$.
 Thus $(b,m)=(mC,m)=f(m)$.
 Similarly, pick $(w,e) \in p^{-1}(F_{m_1})$.
 Then since $(w,e) \in P_1$, $r(w)=eC$. Recall that $r(w)$ is defined as $\sigma_1(w)C$.
 From $eC=\sigma_1(w)C$, there exists $c \in C$ such that $e=c \sigma_1(w)$.
 Therefore, $(w,e)=(w,c \sigma_1(w))=g\left((c,w)\right)$.

 Note that for each $c \in C$, $f(c)=(C,c)$ and $g\left((c,1)\right)=(1,c)$.
 Since the first factor of $P_1$ is $\Gamma = B * F_{m_1}$ and $B$ is defined as $B=H/C$, we can identify $f(c)=g\left((c,1)\right)=(1,c)$ for all $c \in C$.
 Using $f$ and $g$, we obtain a homomorphism 
 $$\Theta: H*_C(C \times F_{m_1} ) \to P_1.$$
 We will verify that $\Theta$ is bijective.
 \begin{enumerate}
    \item $\Theta$ is injective; \\
    Recall that $C=Z(H)$. Since $C$ is central in both factors, $C$ is central in the amalgam $H*_C(C \times F_{m_1} )$.
    In particular, $C$ is a normal subgroup of $H*_C(C \times F_{m_1} )$.
    The quotient group is given by
    $$ H*_C(C \times F_{m_1} ) / C \cong (H/C)* \left ((C \times F_{m_1})/C \right)=B*F_{m_1}=\Gamma,$$
    where the first isomorphism follows from the group presentation.
    For $h \in H$, $p(\Theta(h))=p(hC,h)=hC$ (recall that $p:P_1 \to B*F_{m_1}$ is the projection onto the first factor).
    For $(c,w) \in C \times F_{m_1}$, $p\left(\Theta\left( (c,w) \right)\right)=p(w,c \sigma_1(w))=w$.
    The above calculation implies that the quotient homomorphism $H*_C(C \times F_{m_1} ) \to \Gamma$ is the composition $p \circ \Theta$.
    Therefore, the kernel of $p \circ \Theta$ is $C$.
    So, if $\Theta(v)=1$, then $v \in C \subset H*_C(C \times F_{m_1} )$.
    However, for $c \in C$, $\Theta(c)=(1,c)$ so $\Theta|_C$ is injective.
    Thus, $v=1$, which means that $\Theta$ is injective.
    \item $\Theta$ is surjective; \\
    Let $J:=\Theta(H*_C(C \times F_{m_1} ))$. We want to show that $J=P_1$.
    Clearly, $J \subset P_1$. 
    We will prove the opposite inclusion.
    Note that $p(J)$ contains both $B$ and $F_{m_1}$ so $p(J) = \Gamma$.
    By the definition of $P_1$ and $p$, we have $\ker(p)=\{ (1,c) : c \in C \}.$
    Recall that $ (1,c)=g\left( (c,1) \right)=\Theta \left( (c,1) \right)$, so 
    $\ker(p) \subset J$.
    Then for arbitrary $a \in P_1$, choose $w \in J$ with $p(w)=p(a)$.
    Then $w^{-1}a \in \ker(p) \subset J$, thus $a \in J$.
 \end{enumerate}
 Therefore, $P_1$ is isomorphic to $H*_C(C \times F_{m_1} )$.
 \end{proof}

 Since both $H$ and $F_{m_1}$ are finitely generated (see Lemma~\ref{lem:H-fgrp}), we obtain the following.

\begin{prop} \label{prop:A f.g.}
 The group $P_1$ is finitely generated.
\end{prop}
\begin{rmk}
Lemma~\ref{lem:P-f.g.} cannot be applied directly to prove that $P_1$ is finitely generated, since $A=M/C$ is not finitely presented.
 If $A$ is finitely presented, then $C$ is finitely normally generated in $M$.
 Since $C$ lies in the center $Z(M)$, this means that $C$ is finitely generated.
 However, $C \cong \oplus_\ZZ \ZZ / 2 \ZZ$ is not finitely generated.
\end{rmk}

Due to Proposition~\ref{Fiber_profinite_completion}, we then obtain the isomorphism 
$$ \widehat{P_2} \cong \widehat{P_1 \times F_{m_2}}. $$
We shall show that these two groups $P_2$ and $P_1 \times F_{m_2}$ are the desired counterexample.

\subsection{Being finitely generated and residually finite}

 Now we will verify that $P_2$ and $P_1 \times F_{m_2}$ are finitely generated and residually finite.
 By Proposition~\ref{prop:A f.g.}, $P_1 \times F_{m_2}$ is finitely generated.
 The residual finiteness of $P_1 \times F_{m_2}$ follows from Lemma~\ref{lem:A rf}.
 To show this, we need the following results.

\begin{lem} \label{lem:B-linear}
 The group $B$ has a faithful representation into $\GL(3,\CC(t))$.
\end{lem}
\begin{proof}
 Recall that we define $B=H/Z(H)$.
 In Lemma~\ref{lem:H-fgrp}, every element $a \in H$ can be expressed as 
 $a=(h(t),g(t),\overline{f(t)},s^n)$ where $h(t),g(t) \in R$ and $\overline{f(t)} \in R/2R$.
  In particular, Formula~\eqref{eq:operationH} gives us that
  $$Z(H)=\{ (0,0,\overline{f(t)},s^0) \}$$
  where $\overline{f(t)} \in R/2R \cong \oplus_\ZZ \ZZ / 2\ZZ$.
  Thus, $B$ is isomorphic to $\left \{ (h(t),g(t),\overline{0},s^n) \right \}$ where $h(t),g(t) \in R$.
 Notice that the operation on $B$ is given by 
 \begin{align} \label{eq:operationB}
 \begin{split}
&(h_1(t),g_1(t),\overline{0},s^{n_1}).(h_2(t),g_2(t),\overline{0},s^{n_2}) \\
=&(h_1(t)+t^{n_1}h_2(t),g_1(t)+t^{-n_1}g_2(t),\overline{0},s^{n_1+n_2})
  \end{split}
 \end{align}
 
  From this expression, one can easily verify that $B \cong L/Z(L)$.
  Define $$\Psi:L \to \GL(3,\CC(t)), \quad (h(t),g(t),f(t),s^n) \mapsto \begin{bmatrix}
  t^n & 0 & h(t) \\ 0 & t^{-n} & g(t) \\ 0 & 0 & 1
  \end{bmatrix}.$$
  Then $\Psi$ is a group homomorphism.
  This fact follows from Formula~\eqref{eq:operationL} and matrix multiplication.
  Since $\Psi\left( (h(t),g(t),f(t),s^n) \right)=I_3$ if and only if $g(t)=h(t)=n=0$, the kernel of $\Psi$ is given by $\{ (0,0,f(t),s^0) \}$.
  Since $B$ is isomorphic to $\left \{ (h(t),g(t),\overline{0},s^n) \right \} \cong L/\ker(\Psi)$, we have an injective homomorphism from $B$ into $\GL(3,\CC(t))$.
\end{proof}

 \begin{lem} \label{lem:H into 2copy}
 Let $\overline{L}$ be the image of $L$ under the reduction
 $$ L < \GL(3,R) \twoheadrightarrow \GL(3,\mathbb{F}_2[t,t^{-1}]) < \GL(3,\mathbb{F}_2(t)). $$
 Then there exists an injective homomorphism $\rho: H \to B \times \overline{L}$.
 \end{lem}
 \begin{proof}
 Note that we have two homomorphisms $f_1:L \to B \cong L/Z(L)$ and $f_2:L \to \overline{L}$.
 Using these two homomorphisms, we obtain $j:L \to B \times \overline{L}$.
 Since $\ker(j)=\ker(f_1) \cap \ker(f_2)$ and $\ker(f_1)=Z(L)$, we obtain $\ker(j)=2Z(L)$.
 From $H=L/2Z(L)$, we get an injective homomorphism $\rho: H \to B \times \overline{L}$.
 \end{proof}
 
 Now we prove that $P_1$ is residually finite.

 \begin{lem} \label{lem:A rf}
 The group $P_1$ is residually finite.
 \end{lem}
 \begin{proof}
 Recall the following amalgam expression of $P_1$ (see Lemma~\ref{lem:A amalgam}),
 $$ P_1 \cong H*_C(C \times F_{m_1} ).$$
 Note that the identity map $id_H:H \to H$ and the projection map $C \times F_{m_1} \to C$ agree on $C<H$.
 They induce a retraction $f:P_1 \to H$.
 Combining $f$ with the projection $p:P_1 \to \Gamma$, we obtain a homomorphism $$g:P_1 \to H \times \Gamma.$$
 Note that $\ker(g)=\ker(f) \cap \ker(p)$.
 Recall that $\ker(p)=\{ (1,c): c \in C \}$ but $f|C$ is injective since $C \subset H$.
 Hence, $g$ is injective.
 This means that we can consider $P_1$ as a subgroup of $H \times \Gamma$.

 By Lemma~\ref{lem:H into 2copy}, $H$ is a subgroup of $B \times \overline{L}$.
 Clearly, $\overline{L}$ is finitely generated linear, so $\overline{L}$ is residually finite.
 From Lemma~\ref{lem:B-linear}, we can conclude that $B$ is also finitely generated linear.
 Thus, $B$ is residually finite. This implies that $H<B \times \overline{L}$ is residually finite.

 Recall that $\Gamma=B * F_{m_1}$.
 Since $B$ has a faithful representation into $\GL(3,\CC(t))$ (Lemma~\ref{lem:B-linear}) and $B$ is finitely generated, $B$ is residually finite.
 Since $\Gamma$ is the free product of two residually finite groups, $\Gamma$ is also residually finite \cite[Corollary (ii) on p.44]{MR87652}.
 Therefore, $P_1 < H \times \Gamma$ is residually finite.
 \end{proof}
 
 Therefore, two groups $P_2 < P_1 \times F_{m_2}$ are also residually finite.
 The only remaining thing is the following result.

 \begin{lem}
 The fiber product $P_2$ is finitely generated.
 \end{lem}
 \begin{proof}
 Recall that $P_2$ is the fiber product of $q: P_1 \to M$ and $\sigma_2:F_{m_2} \to M$.
 Also, recall that $P_1$ is finitely generated (Proposition~\ref{prop:A f.g.}) and $M=V(K)$ is finitely presented (Proposition~\ref{prop:G-V(G)}).
 From Lemma~\ref{lem:P-f.g.}, $P_2$ is finitely generated.
 \end{proof}

\subsection{Representations of two groups}

 In this subsection, we will verify that $P_2$ admits a faithful representation into $\GL(d,\CC)$, whereas $P_1 \times F_{m_2}$ is not linear.

 \begin{prop} \label{prop:faithful Sec3}
 The fiber product $P_2$ admits a faithful representation into $\GL(d,\CC)$.
 \end{prop}
 \begin{proof}
 Since $P_2$ is the fiber product of $q: P_1 \to M$ and $\sigma_2:F_{m_2} \to M$, $$P_2=\left \{ (a,f) \in P_1 \times F_{m_2} : q(a)=\sigma_2(f) \right \}. $$
 Also, $P_1$ is the fiber product of $r:\Gamma \to A$ and $M \to A=M/C$.
 Therefore, we obtain
 \begin{align*}
 P_2 &= \left \{ (a,f) \in P_1 \times F_{m_2} : q(a)=\sigma_2(f) \right \} \\
 &\cong \left \{ (\gamma,e,f) \in \Gamma \times M \times F_{m_2} : q((\gamma,e))=\sigma_2(f), r(\gamma)=eC \right \} \\
 &\cong \left \{ (\gamma,e,f) \in \Gamma \times M \times F_{m_2} : e=\sigma_2(f), r(\gamma)=eC \right \}  \\
 &\cong \left \{ (\gamma,f) \in \Gamma \times F_{m_2}:r(\gamma)=\sigma_2(f)C  \right \}.
 \end{align*}
 The above implies that $P_2$ is a subgroup of $\Gamma \times F_{m_2} = (B*F_{m_1}) \times F_{m_2}$.
 By Lemma~\ref{lem:B-linear}, $B$ has a faithful representation into $\GL(3,\CC(t))$.
 Since a free group admits a faithful representation into $\GL(2,\CC(t))$,
 $\Gamma=B*F_{m_1}$ has a faithful representation into $\GL(d',k)$ for some $d'>0$ and some field $k$ with characteristic $0$ \cite{minty1997faithful}.
 Hence, $\Gamma \times F_{m_2}$ admits a faithful representation into $\GL(d'',k)$ for some $d''>0$ and some field $k$ with characteristic $0$.
 Since $\Gamma \times F_{m_2}$ is finitely generated, $\Gamma \times F_{m_2}$ has a faithful representation into $\GL(d,\CC)$ for some $d>0$.
 Therefore, $P_2$ also admits a faithful representation into $\GL(d,\CC)$ for some $d>0$.
 \end{proof}

 Now we complete the proof of the main result, Theorem~\ref{thm:main1}.

 \begin{prop} \label{prop:Not_linear}
 The group $P_1 \times F_{m_2}$ is not linear.
 \end{prop}
 \begin{proof}
 In order to prove the proposition, it suffices to show that $H$ is not linear since $H<P_1 \cong H*_C(C \times F_{m_1})$ (Lemma~\ref{lem:A amalgam}).

 By construction, $H$ contains $C=Z(H)$ and $C$ is isomorphic to $R/2R \cong \oplus_\ZZ \ZZ/2\ZZ$.
 However, if $G < \GL(d,F)$ has finite exponent where the characteristic of $F$ is $0$, then $G$ must be finite \cite[p.15]{wehrfritz2012infinite}.
 Thus, $H$ does not admit a faithful representation into $\GL(d,F)$ with $\textnormal{char}(F)=0$.

 Next, assume that $\textnormal{char}(F)=p \geq 3$ and there exists a faithful representation $\rho:H \to \GL(d,F)$.
 We may assume that $F$ is algebraically closed.
 For $c \in C$, $\rho(c)^2=I$.
 Since $p \neq 2$, the minimal polynomial of $\rho(c)$ divides $X^2-1$, and $X^2-1=(X+1)(X-1)$ has two distinct roots.
 Thus, $\rho(c)$ is diagonalizable with the eigenvalues $\pm 1$.
 Since $C \cong \oplus_\ZZ \ZZ/2\ZZ$, choose a subgroup $C_{d+1} \cong (\ZZ / 2\ZZ)^{d+1}$.
 Since $\rho(C_{d+1})$ is abelian, $\rho(C_{d+1})$ can be simultaneously diagonalized.
 Thus, up to conjugation, $\rho(C_{d+1}) \subset \left \{ \diag(\epsilon_1 , \cdots, \epsilon_d) : \epsilon_i = \pm1 \right \}.$
 However, the cardinalities of $C_{d+1}$ and $\left \{ \diag(\epsilon_1 , \cdots, \epsilon_d) : \epsilon_i = \pm1 \right \}$ are $2^{d+1}$ and $2^d$, respectively.
 Hence, $\rho$ is not injective, and $H$ does not have a faithful representation into $\GL(d,F)$ with odd characteristic.

 Lastly, assume that $\textnormal{char}(F)=2$ and $H$ has a faithful representation $\rho:H \to \GL(d,F)$.
 In this part, we use the properties of B established in Lemma~\ref{lem:PropertyofB}.
 First, note that $H$ is solvable since $H$ is a central extension of $B$ and $B$ is solvable (due to Lemma~\ref{lem:PropertyofB}).
 By Theorem 4.4 in \cite{detinko2019linear}, there exists a finite index subgroup $H_0$ of $H$ such that $\rho(H_0)$ is unipotent-by-abelian.
 Note that any element in $\rho([H_0,H_0])$ is unipotent, so for each $U \in \rho([H_0,H_0])$, we can express it as $U=I+N$ where $N^d=O$.
 For $k \in \NN$ with $2^k \geq d$, we get $U^{2^k} = (I+N)^{2^k}=I+N^{2^k}=I$.
 This implies that $\rho([H_0,H_0])$ is a torsion group, and so is $[H_0,H_0]$.
 Now consider the surjection $\pi: H \to B \cong H/Z(H)$, and let $B_0=\pi(H_0)$.
 Then $B_0$ is a finite index subgroup of $B$.
 Since $[H_0,H_0]$ is a torsion group, so is $[B_0,B_0]$.
 However, $B$ is torsion-free by Lemma~\ref{lem:PropertyofB}, so $[B_0,B_0]=1$.
 This means that $B_0$ is abelian, thus $B$ is virtually abelian.
 This contradicts Lemma~\ref{lem:PropertyofB}.
 Therefore, $H$ does not admit a faithful representation $\rho:H \to \GL(d,F)$ with characteristic $2$.
 \end{proof}

 \begin{lem} \label{lem:PropertyofB}
 The group $B$ is torsion-free solvable and not virtually abelian.
 \end{lem}
 \begin{proof}
 Recall that $B=H/Z(H)$ and $B$ is isomorphic to $\left \{ (h(t),g(t),\overline{0},s^n) \right \}$ where $h(t),g(t) \in R$ with the operation given by Formula~\eqref{eq:operationB}.
 From this, one can check that $$B \cong (R \times R) \rtimes \ZZ,$$
 so $B$ is solvable.
 Since $R$ and $\ZZ$ are torsion-free, so is $B$.

 Now, choose an arbitrary finite index subgroup $B_0<B$.
 Then for some positive integers $m,n$, we have $(m,0,\overline{0},s^0), (0,0,\overline{0},s^n) \in B_0$.
 However,
 $$ (m,0,\overline{0},s^0)(0,0,\overline{0},s^n) = (m,0,\overline{0},s^n)
 \neq (t^nm,0,\overline{0},s^n) = (0,0,\overline{0},s^n)(m,0,\overline{0},s^0), $$
 so $B_0$ cannot be abelian.
 Since $B_0$ is arbitrary, $B$ is not virtually abelian.
 \end{proof}

 By putting $G_1:=P_2$ and $G_2:=P_1 \times F_{m_2}$, we establish Theorem~\ref{thm:main1}.

 \begin{rmk}
 In Lemma~\ref{lem:H into 2copy}, we embed $H$ into the direct product of two linear groups, but the base fields have different characteristics.
 \end{rmk}

 \begin{rmk}
 The failure of faithfulness can be seen explicitly.
 Let $p: P_1 \to \Gamma$ be the projection onto the first factor.
 Choose a faithful representation $\sigma:\Gamma \times F_{m_2} \to \GL(d,\CC)$ as in the proof of Proposition~\ref{prop:faithful Sec3}.
 Then $\widetilde{\sigma}:=\sigma \circ (p \times id): P_1 \times F_{m_2} \to \GL(d,\CC)$ restricts to a faithful representation of $P_2$.
 Note that  
 $$\ker(\widetilde{\sigma}) \cong C \cong \oplus_\ZZ \ZZ / 2\ZZ,$$
 since $\ker(p) \cong C$.
 By Theorem~\ref{thm:Grothendieck}, $\widetilde{\sigma}$ is the unique extension of the restriction $\widetilde{\sigma}|_{P_2}$.
 \end{rmk}

\section{Faithful representations over fields of positive characteristic} \label{sec:main2}

 We dedicate this section to discussing faithful representations into $\GL(d,k)$ with positive characteristic.
 More precisely, our aim in this section is to prove the following result.
 
\begin{thm} \label{thm:main2}
There exist two finitely generated residually finite groups $H_1$ and $H_2$ such that:
\begin{enumerate}[label=(\arabic*)]
  \item $\widehat{H_1} \cong \widehat{H_2}$;
  \item For every field $k$ with characteristic $p \geq 5$ and transcendence degree $r \geq 2$,
there exists $d>0$ such that $H_1$ has a faithful representation into $\GL(d,k)$;
  \item $H_2$ cannot admit a faithful representation into $\GL(d,L)$ for any $d>0$ and any field $L$ with characteristic $p > 0$.
\end{enumerate}
\end{thm}

 We point out that the construction in Theorem~\ref{thm:main2} is entirely independent of that in Theorem~\ref{thm:main1}.
 Therefore, we will reuse some notation from Section~\ref{sec:main1}.

\subsection{The construction of two groups with isomorphic profinite completions}

 We will construct a fiber product with the desired property.
 For this, we start with the universal central extension of Thompson's group $T$ introduced in Section~\ref{subsec:UCET}.
 Let $T$ be Thompson's group, and consider the universal central extension
 $$ 1 \to \ZZ^2 \to E \overset{\pi}{\to} T \to 1.$$
 Let $\{z, w\}$ be a basis of $\ZZ^2$.

\begin{lem} \label{Lem:Fiber_product1}
 There exists a surjective homomorphism $\textnormal{pr} : \pi_1(S_g) \to T$ for some $g$, satisfying $$(1,z) \in Z(A) \cap [A,A]$$
 where $A$ is the fiber product of $\textnormal{pr}:\pi_1(S_g) \to T$ and $\pi:E \to T$.
 Here, $S_g$ is a closed orientable surface of genus $g$.
\end{lem}
\begin{proof}
 Since $E$ is perfect, we have
 $$ z = \prod_{i=1}^m [x_i,y_i] $$
 where $x_i,y_i \in E$.
 Since $T$ is finitely generated, let $\{ t_1, \cdots, t_s \}$ be a finite generating set for $T$.
 Let $g:=m+s$, and define $\textnormal{pr} : \pi_1(S_g) \to T$ by
\begin{align} \label{Form-Surj.}
\begin{split}
a_i \mapsto \pi(x_i), \quad &b_i \mapsto \pi(y_i) \qquad (1 \leq i \leq m) \\
a_{m+j} \mapsto t_j, \quad &b_{m+j} \mapsto 1 \qquad (1 \leq j \leq s)
\end{split}
\end{align}
where $\{a_i, b_i : 1 \leq i \leq g \}$ is the standard generating set for $\pi_1(S_g)$.
  Note that $\textnormal{pr}$ is well-defined and surjective, so we can define the fiber product 
 $$ A:= \left\{ (\gamma,e) : \textnormal{pr}(\gamma)=\pi(e) \right\}<\pi_1(S_g) \times E. $$

 The only remaining thing to prove is $(1,z) \in Z(A) \cap [A,A]$.
 From Formula~\eqref{Form-Surj.}, $(a_i,x_i), (b_i,y_i) \in A$ for $1 \leq i \leq m$, and $(a_{m+j},\tilde{t_j}), (b_{m+j},1) \in A$ for $1 \leq j \leq s$. 
 Here, $\tilde{t_j} \in E$ with $\pi(\tilde{t_j})=t_j$.
 By letting $x_{m+j}=\tilde{t_j}$ and $y_{m+j}=1$, we obtain
 \begin{align} \label{commutator-z}
  \prod_{i=1}^g \left[ (a_i,x_i),(b_i,y_i) \right]=(1,z).
  \end{align}
 This proves $(1,z) \in [A,A]$.
 Since $z \in Z(E)$, clearly $(1,z) \in Z(A)$.
\end{proof}

 Let $q_1:A \to E$ be the projection of $A$ onto the second factor $E$.
 Since $E$ is finitely presented (Proposition~\ref{prop:E-property}), there exists a surjective homomorphism $q_2:F_n \to E$ for some $n$.
 Now we consider the fiber product $P$ of $q_1, q_2$, namely,
 $$P:=\left \{ (a,u) : q_1(a)=q_2(u) \right \}< A \times F_n.$$
 To apply Proposition~\ref{Fiber_profinite_completion}, we need the following.

 \begin{lem} \label{lem:A-f.g. 4}
 The fiber product $A$ is finitely generated.
 \end{lem}
 \begin{proof}
 Recall that $A$ is the fiber product defined by $$ A:= \left\{ (\gamma,e) : \textnormal{pr}(\gamma)=\pi(e) \right\}<\pi_1(S_g) \times E. $$
 Using the projection of $A$ onto the first factor $\pi_1(S_g)$, we obtain the following short exact sequence
 $$ 1 \to \ker(E \to T) \to A \to \pi_1(S_g) \to 1. $$
 Recall that $\ker(E \to T) \cong \ZZ^2$.
 Since both $\ZZ^2$ and $\pi_1(S_g)$ are finitely generated, so is $A$.
\end{proof}
 
Therefore, Proposition~\ref{prop:E-property}, Lemma~\ref{lem:E:trivial_completion}, and Proposition~\ref{Fiber_profinite_completion} give us an isomorphism between the profinite completions, $$\widehat{P} \cong \widehat{A \times F_n}.$$

\subsection{Being finitely generated and residually finite}

 We then check that the two groups $P$ and $A \times F_n$ are finitely generated and residually finite.
First, from Lemma~\ref{lem:A-f.g. 4}, $A \times F_n$ is finitely generated.
We shall show that $P$ is finitely generated.

\begin{lem}
$P$ is finitely generated.
\end{lem}
\begin{proof}
We can describe $P$ as 
\begin{align*}
 P &= \left \{ (a,u) \in A \times F_n: q_1(a)=q_2(u) \right \} \\
 &\cong \left \{ (\gamma,e,u) \in \pi_1(S_g) \times E \times F_n : \textnormal{pr}(\gamma)=\pi(e), q_1((\gamma,e))=q_2(u) \right \} \\
 &\cong \left \{ (\gamma,e,u) \in \pi_1(S_g) \times E \times F_n : \textnormal{pr}(\gamma)=\pi(e), e=q_2(u) \right \} \\
 &\cong \{ (\gamma,u) \in \pi_1(S_g) \times F_n : \textnormal{pr}(\gamma)=\pi(q_2(u)) \}. 
\end{align*}

Thus,  $P$ is indeed the fiber product of $\textnormal{pr}:\pi_1(S_g) \to T$ and $\pi \circ q_2  : F_n \to T.$
Since $\pi_1(S_g)$ and $F_n$ are finitely generated and $T$ is finitely presented, $P$ is finitely generated due to Lemma~\ref{lem:P-f.g.}.
\end{proof}

Therefore, the only remaining thing to prove is that $A \times F_n$ is residually finite.
Clearly, $F_n$ is residually finite.
Since every finitely generated linear group is residually finite, we prove the following.

\begin{prop} \label{A-R representation}
 $A$ has a faithful representation into $\GL(d,\RR)$.
\end{prop}
\begin{proof}
  We claim that $A$ has the following group presentation
 $$ B := \left< \alpha_1,\beta_1, \cdots , \alpha_g,\beta_g,z_1,z_2 : z_1,z_2 \textnormal{ central, } \prod_{i=1}^g [\alpha_i,\beta_i]=z_1  \right >. $$
 Define $ \phi:B \to A $ by
 $$\alpha_i \mapsto (a_i,x_i), \quad \beta_i \mapsto (b_i,y_i), \quad z_1 \mapsto (1,z), \quad z_2 \mapsto (1,w).$$
 Then, clearly, $(1,z),(1,w)$ are central in $A$. 
 By Equation~\eqref{commutator-z}, we have
 $$ \prod_{i=1}^g [\phi(\alpha_i),\phi(\beta_i)] = (1,z) = \phi(z_1). $$
 Thus, $\phi: B \to A$ is a well-defined homomorphism.
 Moreover, $\phi$ is surjective.
 From the short exact sequence $$ 1 \to \ZZ^2 \to A \to \pi_1(S_g) \to 1$$ in Lemma~\ref{lem:A-f.g. 4},
 we conclude that $(a_i,x_i), (b_i,y_i), (1,z), (1,w)$ generate $A$.
 We want to check that $\phi$ is injective.
 Note that $B/\left< z_1,z_2\right> \cong \pi_1(S_g)$ and 
 $A/\left< \phi(z_1),\phi(z_2) \right> \cong \pi_1(S_g)$.
 Since $\pi_1(S_g)$ is Hopfian, we obtain $\ker(\phi) \subset \left< z_1, z_2 \right>$.
 However, $\phi(z_1^m z_2^n)=(1,z^mw^n)$.
 Since $\{z,w\}$ is a basis for $\ZZ^2$, $\ker(\phi)$ is trivial, so $\phi$ is an isomorphism.

 From the presentation, we obtain 
 $$A \cong B \cong C \times \ZZ $$ where $C = \left< \alpha_1,\beta_1, \cdots , \alpha_g,\beta_g,z_1 : z_1 \textnormal{ central, } \prod_{i=1}^g [\alpha_i,\beta_i]=z_1  \right >.$
 By Lemma~4.2 in \cite{ye2025linearity}, there exists an injective homomorphism $i:C \to \GL(d',\RR)$.
 Thus, we have an injective homomorphism $\tilde{i}: A \to \GL(d,\RR)$ for some $d>0$.
\end{proof}

\begin{rmk}
 Lemma~4.2 of \cite{ye2025linearity} states only that $C$ is linear, but its proof yields a faithful representation into $\GL(d,\RR)$.
 To prove Lemma~4.2 in \cite{ye2025linearity}, they used the fact that $C$ can be embedded into the direct product of a surface group and a finitely generated nilpotent group.
 Both groups admit a faithful representation into $\GL(d,\RR)$ for some $d>0$.
\end{rmk}

\subsection{Representations of two groups}

Let $k$ be a field with characteristic $p \geq 5$ and transcendence degree $r \geq 2$.
We will show that $P$ has a faithful representation into $\GL(d,k)$ but $A \times F_n$ does not have a faithful representation into $\GL(d,L)$ for any field $L$ with positive characteristic.

\begin{prop}
 The fiber product $P$ has a faithful representation into $\GL(d,k)$.
\end{prop}
\begin{proof}
Recall that $P$ is isomorphic to
$$ P \cong \{ (\gamma,u) \in \pi_1(S_g) \times F_n : \textnormal{pr}(\gamma)=\pi(q_2(u)) \}.$$
By Theorem 5 in \cite{fisher2018character}, there exists an injective homomorphism $i:\pi_1(S_g) \to \GL(d',k)$ for some $d'>0$.
Since $F_n < \pi_1(S_g)$, we have an injective homomorphism $P \to \GL(d,k)$ for some $d>0$.
\end{proof}

 Since $k$ is arbitrary, we conclude that for any field $k$ with characteristic $p \geq 5$ and transcendence degree $r \geq 2$, the group $P$ has a faithful representation into $\GL(d,k)$ for some $d>0$.

\begin{rmk}
 To apply Theorem 5 in \cite{fisher2018character}, $g \geq 2$ is required.
 Based on our construction, $g \geq 2$ follows immediately.
 Recall that $\textnormal{pr} : \pi_1(S_g) \to T$ is surjective, so the cases of $g=0$ and $g=1$ are impossible.
\end{rmk}

 Now, we establish an obstruction to being linear for $A \times F_n$.
 We note that this obstruction is valid for any positive characteristic.
 We begin with the following statement.

\begin{lem}[Theorem 3 in \cite{button2016mapping}] \label{Lem:field_p-1}
 Let $L$ be a field with $\textnormal{char}(L)>0$, and let $G < \GL(d,L)$.
 Choose an element $ g \in G$ of infinite order.
 Then the image of $g$ in the abelianization of the centralizer $C_G(g)$ also has infinite order.
\end{lem}

 We use the above lemma to establish a criterion for non-linearity.

\begin{lem} \label{Lem:Obstruction1}
 Let $L$ be a field with $\textnormal{char}(L)>0$.
 For $H < \GL(d,L)$, every element in $Z(H) \cap [H,H]$ is a torsion element.
\end{lem}
\begin{proof}
 Suppose that $a \in Z(H) \cap [H,H]$ has infinite order. Since $a \in Z(H)$, the centralizer $C_H(a)$ is $H$.
 By Lemma~\ref{Lem:field_p-1}, $a$ has infinite order in $H/[H,H]$.
 However, $a \in [H,H]$, so the image of $a$ in $H/[H,H]$ is trivial.
\end{proof}

 Now we complete the proof of Theorem~\ref{thm:main2}, the main theorem in this section.

\begin{prop}
 $A \times F_n$ does not have a faithful representation into $\GL(d,L)$ for any field $L$ with $\textnormal{char}(L)>0$.
\end{prop}
\begin{proof}
 It suffices to show that $A$ does not admit a faithful representation into $\GL(d,L)$.
 By Lemma~\ref{Lem:Fiber_product1}, $(1,z) \in Z(A) \cap [A,A]$.
 Since $z \in \ZZ^2$ is of infinite order, $(1,z)$ has infinite order.
 However, by Lemma~\ref{Lem:Obstruction1}, when $A$ has a faithful representation into $\GL(d,L)$,
 the subgroup $Z(A) \cap [A,A]$ cannot contain an element of infinite order.
\end{proof}

 By letting $H_1:=P$ and $H_2:=A \times F_n$, we establish Theorem~\ref{thm:main2}.

\begin{rmk}
 Note that both $P$ and $A \times F_n$ have a faithful representation into $\GL(d,\RR)$.
 By Proposition~\ref{A-R representation}, $A$ has a faithful representation into $\GL(d,\RR)$ and so does $A \times F_n$.
 Since $P< A \times F_n$, $P$ also has a faithful representation into $\GL(d,\RR)$.
 Thus, the groups $P$ and $A \times F_n$ cannot serve as a counterexample to Question~\ref{conj:Lpp}.
\end{rmk}



\bibliography{preprint}

@book {MR2599132,
    AUTHOR = {Ribes, Luis and Zalesskii, Pavel},
     TITLE = {Profinite groups},
    SERIES = {Ergebnisse der Mathematik und ihrer Grenzgebiete. 3. Folge. A
              Series of Modern Surveys in Mathematics [Results in
              Mathematics and Related Areas. 3rd Series. A Series of Modern
              Surveys in Mathematics]},
    VOLUME = {40},
   EDITION = {Second},
 PUBLISHER = {Springer-Verlag, Berlin},
      YEAR = {2010},
     PAGES = {xvi+464},
      ISBN = {978-3-642-01641-7},
   MRCLASS = {20E18},
  MRNUMBER = {2599132},
       DOI = {10.1007/978-3-642-01642-4},
       URL = {https://doi.org/10.1007/978-3-642-01642-4},
}

@article {MR2119723,
    AUTHOR = {Bridson, Martin R. and Grunewald, Fritz J.},
     TITLE = {Grothendieck's problems concerning profinite completions and
              representations of groups},
   JOURNAL = {Ann. of Math. (2)},
  FJOURNAL = {Annals of Mathematics. Second Series},
    VOLUME = {160},
      YEAR = {2004},
    NUMBER = {1},
     PAGES = {359--373},
      ISSN = {0003-486X,1939-8980},
   MRCLASS = {20E18 (20E26)},
  MRNUMBER = {2119723},
MRREVIEWER = {Jan-Christoph\ Schlage-Puchta},
       DOI = {10.4007/annals.2004.160.359},
       URL = {https://doi.org/10.4007/annals.2004.160.359},
}

@article {MR3531665,
    AUTHOR = {Bridson, Martin R.},
     TITLE = {The strong profinite genus of a finitely presented group can
              be infinite},
   JOURNAL = {J. Eur. Math. Soc. (JEMS)},
  FJOURNAL = {Journal of the European Mathematical Society (JEMS)},
    VOLUME = {18},
      YEAR = {2016},
    NUMBER = {9},
     PAGES = {1909--1918},
      ISSN = {1435-9855,1435-9863},
   MRCLASS = {20E18 (20F65 20F67)},
  MRNUMBER = {3531665},
MRREVIEWER = {Martyn\ Quick},
       DOI = {10.4171/JEMS/633},
       URL = {https://doi.org/10.4171/JEMS/633},
}

@article {nikolov2007finitely,
    AUTHOR = {Nikolov, Nikolay and Segal, Dan},
     TITLE = {On finitely generated profinite groups. {I}. {S}trong
              completeness and uniform bounds},
   JOURNAL = {Ann. of Math. (2)},
  FJOURNAL = {Annals of Mathematics. Second Series},
    VOLUME = {165},
      YEAR = {2007},
    NUMBER = {1},
     PAGES = {171--238},
      ISSN = {0003-486X,1939-8980},
   MRCLASS = {20E18 (20E32 20F12)},
  MRNUMBER = {2276769},
MRREVIEWER = {Benjamin\ Klopsch},
       DOI = {10.4007/annals.2007.165.171},
       URL = {https://doi.org/10.4007/annals.2007.165.171},
}

@article{reid2018profinite,
  title={Profinite rigidity},
  author={Reid, Alan W},
  journal={Proc. Int. Cong. of Math.(Rio de Janeiro)},
  volume={2},
  pages={1211--1234},
  year={2018}
}

@article{echtler2024bounded,
  title={Bounded cohomology is not a profinite invariant},
  author={Echtler, Daniel and Kammeyer, Holger},
  journal={Canadian Mathematical Bulletin},
  volume={67},
  number={2},
  pages={379--390},
  year={2024},
  publisher={Canadian Mathematical Society}
}

@article {cheetham2024property,
    AUTHOR = {Cheetham-West, Tamunonye and Lubotzky, Alexander and Reid,
              Alan W. and Spitler, Ryan},
     TITLE = {Property {FA} is not a profinite property},
   JOURNAL = {Groups Geom. Dyn.},
  FJOURNAL = {Groups, Geometry, and Dynamics},
    VOLUME = {19},
      YEAR = {2025},
    NUMBER = {3},
     PAGES = {1081--1087},
      ISSN = {1661-7207,1661-7215},
   MRCLASS = {20E18 (20E06 20E26)},
  MRNUMBER = {4945540},
MRREVIEWER = {Steffen\ Kionke},
       DOI = {10.4171/ggd/802},
       URL = {https://doi.org/10.4171/ggd/802},
}

@inproceedings{platonov1986grothendieck,
  title={On the {G}rothendieck problem of profinite completions of groups},
  author={Platonov, Vladimir Petrovich and Tavgen', Oleg Ignat'evich},
  booktitle={Doklady Akademii Nauk},
  volume={288},
  number={5},
  pages={1054--1058},
  year={1986},
  organization={Russian Academy of Sciences}
}

@article{button2016mapping,
  title={Mapping class groups are not linear in positive characteristic},
  author={Button, JO},
  journal={arXiv preprint arXiv:1610.08464},
  year={2016}
}

@article{heuer2021spectrum,
  title={The spectrum of simplicial volume},
  author={Heuer, Nicolaus and L{\"o}h, Clara},
  journal={Inventiones mathematicae},
  volume={223},
  number={1},
  pages={103--148},
  year={2021},
  publisher={Springer}
}

@article{fisher2018character,
  title={Character varieties and actions on products of trees},
  author={Fisher, David and Larsen, Michael and Spatzier, Ralf and Stover, Matthew},
  journal={Israel Journal of Mathematics},
  volume={225},
  number={2},
  pages={889--907},
  year={2018},
  publisher={Springer}
}

@article{ye2025linearity,
  title={Linearity and virtual poly-freeness of the fundamental group of plane curves of degree at most five},
  author={Ye, Shengkui and Zhu, Kejia},
  journal={arXiv preprint arXiv:2512.08642},
  year={2025}
}

@article{baik2026co,
  title={Co-Hopfianity is not a profinite property},
  author={Baik, Hyungryul and Jang, Wonyong},
  journal={arXiv preprint arXiv:2603.04197},
  year={2026}
}

@article{houcine2007embeddings,
  title={Embeddings in finitely presented groups which preserve the center},
  author={Ould Houcine, Abderezak},
  journal={Journal of Algebra},
  volume={307},
  number={1},
  pages={1--23},
  year={2007},
  publisher={Elsevier}
}

@article{palmer2025embedding,
  title={Embedding groups into acyclic groups},
  author={Palmer, Martin and Wu, Xiaolei},
  journal={arXiv preprint arXiv:2510.16879},
  year={2025}
}

@incollection{thompson1980embeddings,
  title={Embeddings into finitely generated simple groups which preserve the word problem},
  author={Thompson, Richard J},
  booktitle={Studies in Logic and the Foundations of Mathematics},
  volume={95},
  pages={401--441},
  year={1980},
  publisher={Elsevier}
}

@article{wu2025embedding,
  title={Embedding groups into boundedly acyclic groups},
  author={Wu, Fan and Wu, Xiaolei and Zhao, Mengfei and Zhou, Zixiang},
  journal={Journal of the London Mathematical Society},
  volume={111},
  number={5},
  pages={e70164},
  year={2025},
  publisher={Wiley Online Library}
}

@article {rauzy2020computability,
    AUTHOR = {Rauzy, Emmanuel},
     TITLE = {Computability of finite quotients of finitely generated
              groups},
   JOURNAL = {J. Group Theory},
  FJOURNAL = {Journal of Group Theory},
    VOLUME = {25},
      YEAR = {2022},
    NUMBER = {2},
     PAGES = {217--246},
      ISSN = {1433-5883,1435-4446},
   MRCLASS = {20F10 (20E15 20E26 20F05)},
  MRNUMBER = {4388360},
MRREVIEWER = {Denis\ E.\ Serbin},
       DOI = {10.1515/jgth-2020-0029},
       URL = {https://doi.org/10.1515/jgth-2020-0029},
}

@article{minty1997faithful,
  title={Faithful representations of free products},
  author={Minty, R John H},
  journal={Journal of the London Mathematical Society},
  volume={56},
  number={1},
  pages={137--148},
  year={1997},
  publisher={Oxford University Press}
}

@book{wehrfritz2012infinite,
  title={Infinite linear groups: an account of the group-theoretic properties of infinite groups of matrices},
  author={Wehrfritz, Bertram},
  year={2012},
  publisher={Springer Science \& Business Media}
}

@article{detinko2019linear,
  title={Linear groups and computation},
  author={Detinko, Alla S and Flannery, Dane L},
  journal={Expositiones Mathematicae},
  volume={37},
  number={4},
  pages={454--484},
  year={2019},
  publisher={Elsevier}
}

@article{fournier2026profinite,
  title={Profinite isomorphisms, stable commutator length, and fixed point properties},
  author={Fournier-Facio, Francesco},
  journal={arXiv preprint arXiv:2603.12095},
  year={2026}
}

@inproceedings{bridson2026chasing,
  title={Chasing finite shadows of infinite groups through geometry},
  author={Bridson, Martin R},
  booktitle={European Congress of Mathematics},
  pages={1--44},
  year={2026},
  organization={European Mathematical Society-EMS-Publishing House GmbH}
}

@article{sabbagh1991polycyclic,
  title={Polycyclic groups, finite images, and elementary equivalence},
  author={Sabbagh, Gabriel and Wilson, John S},
  journal={Archiv der Mathematik},
  volume={57},
  number={3},
  pages={221--227},
  year={1991},
  publisher={Springer}
}

@article {kionke2023amenability,
    AUTHOR = {Kionke, Steffen and Schesler, Eduard},
     TITLE = {Amenability and profinite completions of finitely generated
              groups},
   JOURNAL = {Groups Geom. Dyn.},
  FJOURNAL = {Groups, Geometry, and Dynamics},
    VOLUME = {17},
      YEAR = {2023},
    NUMBER = {4},
     PAGES = {1235--1258},
      ISSN = {1661-7207,1661-7215},
   MRCLASS = {20E18 (20E08 20E26 43A07)},
  MRNUMBER = {4641368},
MRREVIEWER = {Enrico\ Jabara},
       DOI = {10.4171/ggd/732},
       URL = {https://doi.org/10.4171/ggd/732},
}

@article{aka2012profinite,
  title={Profinite completions and {K}azhdan’s property ({T})},
  author={Aka, Menny},
  journal={Groups, Geometry, and Dynamics},
  volume={6},
  number={2},
  pages={221--229},
  year={2012}
}

@article{cotton2012conjugacy,
  title={Conjugacy separability of 1-acylindrical graphs of free groups},
  author={Cotton-Barratt, Owen and Wilton, Henry},
  journal={Mathematische Zeitschrift},
  volume={272},
  number={3},
  pages={1103--1114},
  year={2012},
  publisher={Springer}
}

@article{lubotzky2014finiteness,
  title={Finiteness properties and profinite completions},
  author={Lubotzky, Alexander},
  journal={Bulletin of the London Mathematical Society},
  volume={46},
  number={1},
  pages={103--110},
  year={2014},
  publisher={Wiley Online Library}
}

@article{bridson2024profinite,
  title={Profinite completions of free-by-free groups contain everything},
  author={Bridson, Martin R},
  journal={The Quarterly Journal of Mathematics},
  volume={75},
  number={1},
  pages={139--142},
  year={2024},
  publisher={Oxford University Press, UK}
}

@article{jang2026bi,
  title={Bi-orderability and generalized torsion elements from the perspective of profinite properties},
  author={Jang, Wonyong and Kim, Junseok},
  journal={Bulletin of the London Mathematical Society},
  volume={58},
  number={3},
  pages={e70308},
  year={2026},
  publisher={Wiley Online Library}
}

@article{jaikin2025some,
  title={Some remarks on {G}rothendieck pairs},
  author={Jaikin-Zapirain, Andrei and Lubotzky, Alexander},
  journal={Groups, Geometry, and Dynamics},
  volume={19},
  number={2},
  pages={597--616},
  year={2025}
}

@article{dixon1982profinite,
  title={Profinite completions and isomorphic finite quotients},
  author={Dixon, John D and Formanek, Edward W and Poland, John C and Ribes, Luis},
  journal={J. Pure Appl. Algebra},
  volume={23},
  number={3},
  pages={227--231},
  year={1982}
}

@article{cannon1996introductory,
  title={Introductory notes on {R}ichard {T}hompson’s groups},
  author={Cannon, James W. and Floyd, W. J. and Parry, W. R.},
  journal={Enseign. Math.(2)},
  volume={42},
  number={3-4},
  pages={215--256},
  year={1996}
}

@article{wu20263,
  title={$3 $-manifold groups are linear},
  author={Wu, Xiaolei and Ye, Shengkui},
  journal={arXiv preprint arXiv:2610.17608},
  year={2026}
}

@article{wise2026linearity,
  title={Linearity of Graph Manifold Groups},
  author={Wise, Daniel T},
  journal={arXiv preprint arXiv:2610.13801},
  year={2026}
}

@article {MR2175160,
    AUTHOR = {Kapovich, Michael},
     TITLE = {Representations of polygons of finite groups},
   JOURNAL = {Geom. Topol.},
  FJOURNAL = {Geometry and Topology},
    VOLUME = {9},
      YEAR = {2005},
     PAGES = {1915--1951},
      ISSN = {1465-3060,1364-0380},
   MRCLASS = {20F69 (20G05)},
  MRNUMBER = {2175160},
MRREVIEWER = {Stephen\ G.\ Brick},
       DOI = {10.2140/gt.2005.9.1915},
       URL = {https://doi.org/10.2140/gt.2005.9.1915},
}

@article{canary2019new,
  title={New nonlinear hyperbolic groups},
  author={Canary, Richard D and Stover, Matthew and Tsouvalas, Konstantinos},
  journal={Bulletin of the London Mathematical Society},
  volume={51},
  number={3},
  pages={547--553},
  year={2019},
  publisher={Wiley Online Library}
}

@article{tholozan2025residually,
  title={Residually finite non-linear hyperbolic groups},
  author={Tholozan, Nicolas and Tsouvalas, Konstantinos},
  journal={Commentarii Mathematici Helvetici},
  volume={100},
  number={1},
  pages={1--9},
  year={2025},
  publisher={European Mathematical Society (EMS)}
}

@article{grothendieck1970representations,
  title={Repr{\'e}sentations lin{\'e}aires et compactification profinie des groupes discrets},
  author={Grothendieck, Alexander},
  journal={Manuscripta mathematica},
  volume={2},
  number={4},
  pages={375--396},
  year={1970},
  publisher={Springer}
}

@book{Weibel_1994,
  place={Cambridge},
  series={Cambridge Studies in Advanced Mathematics},
  title={An Introduction to Homological Algebra}, 
  publisher={Cambridge University Press}, 
  author={Weibel, Charles A.}, 
  year={1994}, 
  collection={Cambridge Studies in Advanced Mathematics}
}

@article {MR87652,
    AUTHOR = {Gruenberg, K. W.},
     TITLE = {Residual properties of infinite soluble groups},
   JOURNAL = {Proc. London Math. Soc. (3)},
  FJOURNAL = {Proceedings of the London Mathematical Society. Third Series},
    VOLUME = {7},
      YEAR = {1957},
     PAGES = {29--62},
      ISSN = {0024-6115,1460-244X},
   MRCLASS = {20.0X},
  MRNUMBER = {87652},
MRREVIEWER = {P.\ Hall},
       DOI = {10.1112/plms/s3-7.1.29},
       URL = {https://doi.org/10.1112/plms/s3-7.1.29},
}

@article{corson2026higman,
  title={Higman--{T}hompson groups and profinite properties of right-angled {C}oxeter groups},
  author={Corson, Samuel M and Hughes, Sam and M{\"o}ller, Philip and Varghese, Olga},
  journal={Selecta Mathematica},
  volume={32},
  number={2},
  pages={28},
  year={2026},
  publisher={Springer}
}
\bibliographystyle{abbrv}

\end{document}